\documentclass[11pt]{article}
\usepackage{amsmath}
\usepackage{amsfonts}
\usepackage{latexsym}
\usepackage{enumerate}
\usepackage{array}
\usepackage{multirow}
\newcommand{\C}{{\bf C}}
\newcommand{\R}{{\bf R}}
\newcommand{\Q}{{\bf Q}}
\newcommand{\Z}{{\bf Z}}

\newcommand{\Aut}{{\rm Aut}}

\newcommand{\rtimes}{\times \kern -2pt
\vrule height 5.2pt depth 0pt width 0.4pt\;}

\begin{document}
\title{Four-dimensional compact solvmanifolds with\\
and without complex analytic structures
}
\author{Keizo Hasegawa\\
Niigata University, JAPAN
}
\date{July, 2003
}

\maketitle
\begin{quote}
Abstract.\,  We classify four-dimensional compact solvmanifolds up to
diffeomorphism, while determining which of them have complex
analytic structures. In particular, we shall see that a four-dimensional
compact solvmanifold $S$ can be written, up to double covering,
as $\Gamma \backslash G$ where $G$ is a simply connected solvable
Lie group and $\Gamma$ is a lattice of $G$, and every complex structure
$J$ on $S$ is the {\em canonical complex structure} induced from a left-invariant
complex structure on $G$. We are thus led to conjecture that complex analytic
structures on compact solvmanifolds are all canonical.
\end{quote}
\medskip

\baselineskip=18pt

\begin{center}
{1.  INTRODUCTION}
\end{center}

In this paper we shall mean by a solvmanifold (nilmanifold) a compact
homogeneous space of solvable (nilpotent) Lie group.
Let $M$ be an $n$-dimensional solvmanifold.
$M$ can be written as $D \backslash G$, where $G$ is a simply connected
solvable Lie group and $D$ is a closed subgroup (which may not be a
discrete subgroup) of $G$. In the case where $M$ is a nilmanifold,
we can assume that $D$ is a discrete subgroup of $G$.
\smallskip

In the author's paper [\ref{H2}] we observed that complex tori and hyperelliptic
surfaces are the only four-dimensional solvmanifolds which admit K\"ahler
structures. In this paper we extend this result to the case of complex
structures; namely we will show the following:
\medskip

\noindent {\bfseries Main Theorem.}
{\em The complex surfaces with diffeomorphism type of solvmanifolds
are all of the complex tori, hyperelliptic surfaces,
Primary Kodaira surfaces, Secondary Kodaira surfaces and Inoue surfaces.}
\medskip

In Section 2 we classify four-dimensional orientable solvmanifolds
up to diffeomorphism, which include all of the orientable
$T^2$-bundles over $T^2$. In Section 3 we study in more details certain
classes of four-dimensional solvmanifolds
which admit complex structures; we show precisely how to construct complex
structures on them. In particular we obtain the classifications of hyperelliptic
surfaces and secondary Kodaira surfaces as solvmanifolds. 
In Section 4 we give the proof of the main theorem. To be more precise, the
proof consists of two steps: First, we show that complex surfaces stated above in
the main theorem are all obtained by constructing complex structures on certain
four-dimensional solvmanifolds. This is the part we shall see in Section 3.
Next, we show that a complex surface with diffeomorphism type of solvmanifold
must be one of the surfaces stated in the main theorem. This is
the part we actually show in Section 4. As a consequence of the theorem 
together with the results of Section 3, we can see that complex structures on
a complex surface $S$ with diffeomorphism type of solvmanifold
$\Gamma \backslash G$ (up to finite covering) are exactly those which are canonically
induced from left-invariant complex structures on $G$ (Proposition 4.2), where
$G$ is a simply connected solvable Lie group and $\Gamma$ is a lattice of $G$.
In other word these complex structures can be
defined simply on the Lie algebras of $G$ (see Section 3.8).
In Section 5 (Appendix) we give a complete list of all the complex structures on
four-dimensional compact homogeneous spaces, referring to their
corresponding complex surfaces.
A brief description of the list is provided, though the proof is mostly a
combination of many known results, along with the application of
Kodaira's classification of complex surfaces.
\smallskip

We shall discuss also some other important structural problems on
solvmanifolds: while classifying four-dimensional solvmanifolds
we can see the following:

\begin{list}{}{\topsep=0pt \leftmargin=10pt \itemindent=5pt \parsep=0pt \itemsep=5pt}
\item[1)] Every four-dimensional orientable solvmanifold is real
parallelizable (see Proposition 2.3), which is noted in [\ref{AS1}] without proof.
\item[2)] There exist four-dimensional simply connected, unimodular and
non-nilpotent solvable Lie groups which have no lattices (see Section 2.5).
\item[3)] There exist four-dimensional solvmanifolds which
cannot be written as $\Gamma\backslash G$ where $G$ is a simply
connected solvable Lie group and $\Gamma$ is a lattice of $G$, although
they have the solvmanifolds of the above type as double coverings
(see Example 3.5).

\end{list}
\bigskip

\begin{center}
{2.  FOUR-DIMENSIONAL SOLVMANIFOLDS -- CLASSIFICATION\\
UP TO DIFFEOMORPHISM}

\end{center}

We recall some fundamental results on solvmanifolds.
It is well-known [\ref{M}] that an $n$-dimensional solvmanifold $M$ is a fiber
bundle over a torus with fiber a nilmanifold, which we call the {\em Mostow fibration}
of $M$. In particular we can represent the fundamental group $\Gamma$ of $M$
as a group extension of $\Z^k$ by a torsion-free nilpotent group $N$ of
rank $n-k$, where $1 \le k \le n$ and $k = n$ if and only if $\Gamma$ is ableain:

$$0 \rightarrow N \rightarrow \Gamma \rightarrow \Z^k \rightarrow 0.$$
\noindent Conversely any such a group $\Gamma$ (which is called the
{\em Wang group}) can be the fundamental group of some solvmanifold
([\ref{A}, \ref{WG}]). It is also known [\ref{M}] that two solvmanifolds having
isomorphic fundamental groups are diffeomorphic.
\medskip

From now on we denote by $S$ a four-dimensional orientable solvmanifold,
$\Gamma$ the fundamental group of $S$, and $b_1$ the first Betti number
of $S$. The classification is divided into three cases: 
\smallskip

\begin{list}{}{\topsep=0pt \leftmargin=15pt \itemindent=20pt \parsep=0pt \itemsep=0pt}
\item[(2.1)] $2 \le k \le 4$ and $N$ is abelian
\item[(2.2)] $k = 1$ and $N$ is abelian
\item[(2.3)] $k = 1$ and $N$ is non-abelian
\end{list}
\smallskip

It should be noted that these three cases are not exclusive each other. For instance,
it can be easily seen that $\Gamma$ of  the case $k = 3$ can be expressed as
the case $k = 2$, and also as the case $k = 1$ and $N$ is either abelian or
non-abelian. More precisely, we shall see that the class of solvmanifolds of the
case (2.1) coincides with that of $T^2$-bundles over $T^2$ (see Proposition 2.1),
and most of them belong also to the case (2.2) or (2.3).
\bigskip

\noindent{\large 2.1. \bfseries The case $2 \le k \le 4$ and $N$ is abelian}
\smallskip

In the case $k = 4$, $\Gamma$ is abelian and $S$ is a 4-torus with
$b_1 = 4$. In the case $k = 3$, we have $N=\Z^1$, and since $S$ is orientable the
action $\phi:\Z^3 \rightarrow \Aut(\Z)$ is trivial. Hence $\Gamma$ is a
nilpotent group, which can be canonically extended to the nilpotent Lie group
$G$, and $S= \Gamma \backslash G$ is a nilmanifold with $b_1 = 3$. $S$ has a
complex structure defining a primary Kodaira surface (see Section 3.4). 
In the case $k = 2$, we have $N=\Z^2$, and $S$ is a $T^2$-bundle over $T^2$.

We have the following general result (including non-orientable cases), which appears
to be well-known.
\smallskip

\noindent {\bfseries Proposition 2.1.} {\em The diffeomorphism class of
four-dimensional solvmanifolds with $ 2 \le b_1 \le 4$ coincides with that of
$T^2$-bundles over $T^2$.} 
\smallskip

\noindent {\em Proof.} We know that $\Gamma$ of  the case $k = 3$ can be
expressed as the case $k = 2$, and the corresponding nilmaifold with $b_1 = 3$ 
has a structure of $T^2$-bundle over $T^2$. We also know
[\ref{AS2}] that a solvmanifold of general dimension can has a structure of fiber
bundle over $T^{b_1}$ with fiber a nilmanifold (which may be different from
the Mostow fibration). Hence, the four-dimensional solvmanifolds with
$ 2 \le b_1 \le 4$ are all $T^2$-bundles over $T^2$. Conversely,
for a given $T^2$-bundle $\bar{S}$ over $T^2$ there exits a four-dimensional
solvmanfold $S$ with the same fundamental group as $\bar{S}$.
Since the diffeomorphism type of $T^2$-bundle over $T^2$
is determined only by the fundamental group ([\ref{SF}]), $S$ must be
diffeomorphic to $\bar{S}$.
\hfill $\Box$
\medskip

The classification of $T^2$-bundles over $T^2$ is well-known
([\ref{SF}, \ref{U}]). According to the above proposition, we can classify them as
four-dimensional solvmanifolds. Let us first consider the group extension:
$$0 \rightarrow \Z^2 \rightarrow \Gamma \stackrel{r}{\rightarrow} \Z^2 \rightarrow 0,$$

\noindent where the action $\phi: \Z^2 \rightarrow \Aut(\Z^2)$ is defined by
$\phi(e_1), \phi(e_2), e_1 = (1,0), e_2 = (0,1)$. Since the total space is orientable
$\phi(e_1), \phi(e_2) \in {\rm SL}(2,\Z)$,  and since $\phi(e_1)$ and $ \phi(e_2)$ commute
we can assume that $\phi(e_1) = \pm {\rm I}$ (cf. [\ref{SF}]). 
\medskip

\begin{list}{}{\topsep=0pt \leftmargin=10pt \itemindent=5pt}

\item[ {\bf 1)}] In the case where $\phi(e_1) = {\rm I}$,
by taking a section $s$ of $r$ defined on the first facter of $\Z^2$, we can
construct a subgroup $\Z^3$ of $\Gamma$ with the group extension:
$$0 \rightarrow \Z^3 \rightarrow \Gamma \rightarrow \Z^1 \rightarrow 0.$$
Hence this case is included in the case (2,2).
\item[ {\bf 2)}] In the case where $\phi(e_1) = - {\rm I}$, we can assume that
the characteristic polynomial $\Phi$ of $\phi(e_2)$ has the double root $1$ or
two distinct positive real roots. 
\begin{list}{}{\topsep=0pt \leftmargin=10pt \itemindent=5pt \parsep=0pt \itemsep=5pt}
\item[ {\bf a)}] If $\Phi$ has the double root $1$,  $\phi(e_2)$ is expressed in the
following form:
$$\left(
\begin{array}{cc}
1 & n\\
0 & 1
\end{array}
\right).$$
As in the case where $\phi(e_1) =  {\rm I}$, by taking a section $s$ of $r$,
we can construct a subgroup $N$ of $\Gamma$, which is a torsion-free
nilpotent group with the group extension:
$$0 \rightarrow N \rightarrow \Gamma \rightarrow \Z^1 \rightarrow 0.$$

Hence this case is included in the the case (2,3) (this case is precisely  the case 2b)
of (2,3)).
\item[ {\bf b)}] If $\Phi$ has two distinct positive real roots $a_1, a_2$,
we can similarly construct a subgroup $K$ of $\Gamma$ with the group extension:
$$0 \rightarrow K \rightarrow \Gamma \rightarrow \Z^1 \rightarrow 0,$$
where $K = \Z^2 \rtimes \Z$ (a non-nilpotent solvable group) which can be
extended to the solvable Lie group $G =\R^2 \rtimes \R$ defined by the action
$\psi(t) (x,y) = (e^{t \log a_1} x, e^{t \log a_2} y)$. $\bar{S} = K \backslash G$
is a three-dimensional solvmanifold, and the solvmanifold $S$ corresponding to
$\Gamma$ is a fiber bundle over $T^1$ with fiber $\bar{S}$ where the action
$\phi : \Z \rightarrow {\rm Aut}(\bar{S})$ is canonically induced by
$\phi(m) = (- 1)^m {\rm I}_2 \times {\rm I}_1
\in {\rm Aut}(G)$.
\end{list}
\end{list}
\medskip

\noindent {\large \em Remark 2.2.} The solvmanifolds of the case 2a) and 2b)
above are the examples which cannot be written as $\Gamma \backslash G$,
where $G$ is a simply connected solvable Lie group and $\Gamma$ is
a lattice of $G$. However, it is clear that they have the solvmanifolds of the above
type as double coverings.
\bigskip

\noindent{\large  2.2.  \bfseries The case $k = 1$ and $N$ is abelian}
\smallskip

The group extension is split and thus
determined only by the action $\phi: \Z \rightarrow \Aut(\Z^3)$.
Since $S$ is orientable, $\phi(1) = A \in {\rm SL(3,\Z)}$, and the action $\phi$ is
defined by
$$\phi(1)(e_i) = \sum_{j=1}^{3} a_{ij} e_j, i = 1,2,3,$$
where $e_1 = (1,0,0), e_2 = (0,1,0), e_3 = (0,0,1)$.
We mean by roots of $\phi(1)$ the eigenvalues of $\phi(1)$, the roots of the
characteristic polynomial $\Phi$ of $\phi(1)$. Since $\Phi$ is of the form
$\Phi(t)=t^3 - m t^2 + n t -1$,
where $m, n$ are integers, we can see that $\Phi$ has no double root
except $1$ or $-1$. This is an important Lemma (for which the proof is
quite elementary) though it seems to be unknown.
\medskip

\noindent {\bfseries Lemma 2.2.} {\em Let $\Phi(t)$ be a polynomial of the form
$\Phi(t) = t^3 - m t^2 + n t -1\;(m, n \in \Z)$.
Then it has a real double root $a$ if and only if $a=1$ or $-1$ for which
$\Phi(t)=t^3 -3 t^2 + 3 t -1$ or $\Phi(t)=t^3 + t^2 - t -1$ respectively.}
\smallskip

\noindent {\em Proof}. Assume that $\Phi(t)$ has a double root $a$ and
another root $b$. Then we have that $a^2 b = 1, 2a + b = m, a^2 + 2ab = n$,
from which we deduce that $ma^2 - 2na + 3 = 0, 3a^2 - 2ma + n = 0$;
and thus $2(m^2 - 3n) a = mn - 9$.
If $m^2 = 3n$, then $m = n = 3$ and $a=1$. If $m^2 \not= 3n$, then
we have  that $a = \frac{mn - 9}{2(m^2 - 3n)}$, which is a rational number.
Since we have that $2a + \frac{1}{a^2} =m \in \Z$, $a$ must be $-1$ or $1$.

\hfill $\Box$
\medskip

Now we classify the group extensions and their corresponding solvmanifolds
according to the roots of $\phi(1)$.
\medskip

\begin{list}{}{\topsep=0pt \leftmargin=10pt \itemindent=10pt \parsep=0pt \itemsep=5pt}

\item[ {\bf 1)}] If $\phi(1)$ has three distinct positive real roots $a_1,a_2,a_3$,
we have linearly independent eigen vectors $u_1, u_2, u_3$ of
$a_1, a_2, a_3$ respectively. Let $u_i=(u_{i1}, u_{i2}, u_{i3}),\; i=1,2,3$.
Then  $\{v_1, v_2, v_3 \},\;v_j=(u_{1j}, u_{2j}, u_{3j}), j=1,2,3$ defines an
abelian lattice $\Z^3$ of $\R^3$. We define a solvable Lie group
$G = \R^3 \rtimes \R$,
where the action $\phi: \R \rightarrow \Aut(\R^3)$ is defined by
$$\phi(t)((x, y, z)) = (e^{t \log a_1} x, e^{t \log a_2} y, e^{t \log a_3} z),$$
which is a canonical extension of $\phi$.
Then $\Gamma = \Z^3 \rtimes \Z$ is a lattice of $G$,
and  $S=\Gamma \backslash G$ is a solvmanifold. We can see that
$S$ is a $T^2$-bundle over $T^2$. $S$ has $b_1=2$ for the case
where one of the roots is $1$, and $b_1=1$ for the case where none of the roots
is $1$.

\item[ {\bf 2)}] If $\phi(1)$ has three distinct real roots two of which are
negative, $S$ has a double covering solvmanifold of the
above type 1).
\item[ {\bf 3)}] If $\phi(1)$ has a triple root $1$, taking a suitable basis
$\{u_1,u_2,u_3\}$ of $\R^3$, $\phi(1)$ is expressed in either of the following
forms:
$$\left(
\begin{array}[c]{ccc}
1 & 1 & \frac{1}{2}\\
0 & 1 & 1\\
0 & 0 & 1
\end{array}
\right),\;
\left(
\begin{array}[c]{ccc}
1 & 1 & 0\\
0 & 1 & 0\\
0 & 0 & 1
\end{array}
\right).
$$
Let $G=\R^3 \rtimes \R$, where the action $\bar{\phi}: \R \rightarrow
{\rm Aut}(\R^3)$ is defined by
$$\bar{\phi}(t)=
{\rm exp}\;t
\left(
\begin{array}[c]{ccc}
0 & 1 & 0\\
0 & 0 & 1\\
0 & 0 & 0
\end{array}
\right)=\left(
\begin{array}[c]{ccc}
1 & t & \frac{1}{2} t^2\\
0 & 0 & t\\
0 & 0 & 0
\end{array}
\right)
$$
for the former case, and
$$\bar{\phi}(t)=
{\rm exp}\;t
\left(
\begin{array}[c]{ccc}
0 & 1 & 0\\
0 & 0 & 0\\
0 & 0 & 0
\end{array}
\right)=\left(
\begin{array}[c]{ccc}
1 & t & 0\\
0 & 1 & 0\\
0 & 0 & 1
\end{array}
\right)
$$
for the latter case.
Then, as defined in the case 1), $\{v_0, v_1, v_2, v_3\}$ defines a
lattice $\Gamma$ of $G$, and $S=\Gamma \backslash G$ is a nilmanifold.
We can see that $S$ is a nilmanifold with $b_1=2$ for the former case, and
a nilmanifold with $b_1=3$, a primary Kodaira surface, for the latter case.

\item[ {\bf 4)}] If $\phi(1)$ has a single root $1$, and a double root $-1$ or
non-real complex roots $\beta \,(|\beta|=1)$, $\bar{\beta}$,
then $S$ is a $T^2$-bundle over $T^2$ and has $b_1 = 2$.
In case $\phi(1)$ has a double root $-1$, taking a suitable basis
$\{u_1,u_2,u_3\}$ of $\R^3$, $\phi(1)$ is expressed in either of the following
forms:
$$\left(
\begin{array}[c]{ccc}
-1 & 0 & 0\\
0 & -1 & 0\\
0 & 0 & 1
\end{array}
\right),\;
\left(
\begin{array}[c]{ccc}
-1 & 1 & 0\\
0 & -1 & 0\\
0 & 0 & 1
\end{array}
\right).
$$

Except for the latter case in the above, each corresponding
solvmanifold has a complex structure, defining a hyperelliptic surface.
For the details of hyperelliptic surfaces we refer to Section 3.3.

We can see that the solvmanifold corresponding to the latter case has
a nilmanifold with $b_1=3$ as a double covering. 
\item[ {\bf 5)}] If $\phi(1)$ has a positive real root $a$, and non-real complex
roots $\beta \,(|\beta|=b \not= 1)$, $\bar{\beta}$, then $S$ is
an Inoue surface of type $S$.  For the details we refer to Section 3.6.
\end{list}
\bigskip

\noindent{\large 2.3.  \bfseries The case $k = 1$ and $N$ is non-abelian}

\smallskip

The group extension has a non-abelian kernel $N$ which is a torsion-free nilpotent
group of rank $3$. We can see that the extension is split and is determined only by
the action $\phi: \Z \rightarrow \Aut(N)$.
$\phi(1)$ induces an automorphism $\widetilde{\phi}(1)$ of the center $\Z$
of $N$.  $\phi(1)$ also induce an automorphism
$\widehat{\phi}(1)$ of $\Z^2= N/\Z$,
that is, $\widehat{\phi}(1) \in {\rm GL}(2,\Z)$.
\medskip

\begin{list}{}{\topsep=0pt \leftmargin=10pt \itemindent=5pt}

\item[ {\bf 1)}] In the case where $\widetilde{\phi}(1)= {\rm Id}$ and
${\rm det}\, \widehat{\phi}(1)=1$, we have the following:

\begin{list}{}{\topsep=0pt \leftmargin=10pt \itemindent=5pt \parsep=0pt \itemsep=5pt}

\item[{\bf a)}] If $\widehat{\phi}(1)$ has positive real roots $a \,(\not=1), b$,
then $S$ is an Inoue surface of type $S^+$. We refer to Section 3.7 for the details.

\item[{\bf b)}] If $\widehat{\phi}(1)$ has negative real roots $a\,(\not=-1), b$,
then $S$ has a double covering solvmanifold of the above type a).
\item[{\bf c)}] If $\widehat{\phi}(1)$ has a double root $1$, then  we can assume
that  $\widehat{\phi}(1)$ is of the form
$$\left(
\begin{array}{cc}
1 & n\\
0 & 1
\end{array}
\right).$$
We can see that ${\rm H}^1(\Gamma, \Z)= \Z^2 \oplus \Z / n \Z$.
Hence, $S$ is a $T^2$-bundle over $T^2$, for which $b_1=2$ for
the case $n \not= 0$, and $b_1=3$ ($S$ is a primary Kodaira surface) for the
case $n=0$.

\item[{\bf d)}] If $\widehat{\phi}(1)$ has a double root $-1$ or non-real
complex roots $\alpha \,(|\alpha|=1)$, $\bar{\alpha}$, then $S$ has a double
covering of the type 1c) and has $b_1 = 1$.
In case $\widehat{\phi}(1)$ has a double root $-1$, we can assume that
$\widehat{\phi}(1)$ is of the form
$$\left(
\begin{array}{cc}
-1 & n\\
0 & -1
\end{array}
\right).$$

Except for the case $n \not= 0$ in the above, each corresponding solvmanifold
has a complex structure, defining a secondary Kodaira surface. 
We refer to Section 3.5 for the details of secondary Kodaira surfaces.
\end{list}

\item[ {\bf 2)}] In the case where $\widetilde{\phi}(1)=- {\rm Id}$ and
${\rm det}\, \widehat{\phi}(1)=-1$, it is clear that $S$ has a double covering
of a solvmanifold of type 1). We have actually the following:


\begin{list}{}{\topsep=0pt \leftmargin=10pt \itemindent=5pt \parsep=0pt \itemsep=5pt}
\item[{\bf a)}] If $\widehat{\phi}(1)$ has real roots $a \,(\not=1), b$,
then $S$ is an Inoue surface of type $S^-$ (see Section 3.7).

\item[{\bf b)}] If  $\widehat{\phi}(1)$ has roots $1$ and $-1$, then $S$ is a
$T^2$-bundle over $T^2$ and has $b_1 = 2$.
We refer to Example 3.5 for the details.
\end{list}
\end{list}

\medskip

\noindent {\large \em Remark 2.2.} $\Gamma$ can be expressed both as a group
extension of $\Z^1$ by $\Z^3$ and as that of $\Z^1$ by a non-ableain nilpotent
group $N$ if and only if $\Gamma$ is of the form $N \times \Z^1$.
Correspondingly, $S$ belongs to both of the cases (2,2) and (2,3) if and only if
$S$ is a nilmanifold with $b_1=3$ (a Kodaira surface).
\bigskip

\noindent{\large 2.4. \bfseries Parallelizability of solvmanifolds}
\smallskip

As an application of the classification, we can see the following (cf. [\ref{AS1}]). 
\smallskip

\noindent {\bfseries Proposition 2.3.} {\em Any four-dimensional orientable
compact solvmanifold is real parallelizable.}

\smallskip

\noindent {\em Proof.} In fact, let us consider first a solvmanifold
$\Gamma \backslash G$,
where $G= \R^2 \rtimes \R$ with the action $\phi$ defined by
$$\phi(t) (x,y) = (\cos (\pi t) x - \sin (\pi t) y, \sin (\pi t) x+\cos (\pi t) y),$$
and $\Gamma = \Z^2 \rtimes \Z$, which is a lattice of $G$.
Then we have linearly independent left-invariant vector fields
$X_1, X_2, X_3$ on $G$:
\medskip

$X_1= \cos (\pi t) \frac{\partial}{\partial x} +
\sin (\pi t) \frac{\partial}{\partial y},$
$X_2= -\sin (\pi t) \frac{\partial}{\partial x} +
\cos (\pi t) \frac{\partial}{\partial y},
X_3= \frac{\partial}{\partial t}.$
\medskip

\noindent As we have seen in the above classification, a four-dimensional
solvmanifold $S$ is either of the form $\Gamma \backslash G$, where
$G$ is a simply connected solvable Lie group and $\Gamma$ is a lattice
of $G$, or has a double covering of a solvmanifold $\widehat{S}$ of the above
type. For the former case it is clearly parallelizable, and for the latter case
since the covering is of the same type as the above,  we can construct four
linearly independent vector fields on $\widehat{S}$, which are invariant by the
covering transformation, in the same way as the above case (though
they are not in general left-invariant vector fields on $G$). 
\hfill $\Box$
\medskip

\noindent {\large \em Remark 2.4.}
There exists a five-dimensional solvmanifold which is not real
parallelizable ([\ref{AS1}]).

\bigskip

\noindent {\large  2.5. \bfseries Lattices of solvable Lie groups}.

\smallskip

\begin{list}{}{\topsep=0pt \leftmargin=10pt \itemindent=10pt}
\item[1)]  A simply connected solvable Lie group $G = \R^3 \rtimes \R$,
where the action $\phi: \R \rightarrow \Aut(\R^3)$ is defined by
$\phi(t)((x, y, z)) = (e^{a t} x, e^{a t} y, e^{-2a t} z)$,
$a (\not=0) \in \R$, has no lattices, although it has a geometric structure.
This is due to the Lemma, and the fact that $G$ has a circle
action rotating the first two coordinates. An Inoue surface of type $S$
can be interpreted as a geometry of this type.
(for the details we refer to [\ref{W1}]).
\smallskip

\item[2)] We have another type of a simply connected, unimodular and
non-nilpotent solvable Lie group which has no lattices: Let $\mathfrak g$
be a Lie algebra
$${\mathfrak g} = \{X_1,X_2,X_3,X_4\},$$
where bracket multiplications are defined by
$$ [X_4,X_1]=-2 X_1, [X_4,X_2]=X_2, [X_4,X_3]=X_3,$$
$$ [X_2,X_3]=-X_1, [X_1,X_2]=0, [X_1,X_3]=0.$$
Then it can be seen from our classification that the corresponding Lie
group $G$ has no lattices.

\end{list}
\bigskip

\begin{center}
{3.  FOUR-DIMENSIONAL SOLVMANIFOLDS -- WITH THE\\
CONSTRUCTIONS OF CANONICAL COMPLEX STRUCTURES}

\end{center}
\smallskip

We shall study in this section some important classes of
solvmanifolds in further detail. In particular we study in separate sections
certain classes of solvmanifolds on which we can construct canonical complex
structures, determining some explicit classes of complex surfaces. 
\bigskip

\noindent{\large 3.1. \bfseries Some basic examples} 
\smallskip

\noindent {\bfseries Example 3.1.}
Let $\Gamma=\Z \rtimes \Z$, where the action $\phi:
\Z \rightarrow {\rm Aut}(\Z)$ is defined by $\phi(1)
=-1$. Then $\Gamma$ acts on $\R^2$ as follows:
$$(a,b) \cdot (x,y) = (a+(-1)^b x, b+y),$$
$(x,y) \in \R^2, (a,b) \in \Gamma$. $M=\R^2/\Gamma$
is the Klein Bottle. Let $G=\C \rtimes \R$ be a solvable
Lie group defined as follows:
$$(w,t) \cdot (z,s) = (w+e^{\pi i t} z, t+s).$$

Let $D= \{(p+iu, q) |\, p,q \in \Z, u \in \R\}$ and
$D'=\{(iu,0) |\, u \in \R\}$ be closed
subgroups of $G$. Then we have $D=  D' \rtimes \Gamma$, and
$D \backslash G = (D' \backslash G)/\Gamma
= \R^2/\Gamma$.
\bigskip

\noindent {\bfseries Example 3.2.}
Let $G= \C \rtimes \R$ be a solvable Lie group defined as
follows:
$$(w,t) \cdot (z,s) = (w+e^{2 \pi i t} z ,t+s).$$

$\Gamma = \{(p+qi,r) |\, p, q, r \in \Z\}$ is a uniform
discrete subgroup of $G$, which is abelian. $\Gamma \backslash G$ is
diffeomorphic to $T^3$ (a three-dimensional torus).
\bigskip

\noindent {\bfseries Example 3.3.}
Let $\Gamma= \Z^2 \rtimes \Z$, where the action $\phi:
\Z \rightarrow {\rm Aut}(\Z^2)$ is defined by $\phi(1)= A
\in {\rm SL}(2,\Z)$. Assume that $A$ has complex roots
$\alpha, \bar{\alpha}\, (\alpha \not= \bar{\alpha}),
|\alpha| = 1$ with the eigen vector $(a,b) \in \C^2$ of $\alpha$,
or a double root $-1$ with linearly independent eigen vectors $a, b \in \C$
(considering $\R^2$ as $\C$)
Since the multiplication by $\alpha$ preserves the
lattice spanned by $a$ and $b$, $\alpha$ must be
$e^{i \eta t}, \eta = \frac{2}{3} \pi,
\frac{1}{2} \pi$ or $\frac{1}{3} \pi$. Let $G= \C \rtimes \R$
be a solvable Lie group defined as follows:
$$(w,t) \cdot (z,s) = (w+e^{i \eta t} z, t+s),$$
where $\eta = \pi, \frac{2}{3} \pi, \frac{1}{2} \pi$ or
$\frac{1}{3} \pi$. Considering $\Gamma$ as the lattice of $G$ spanned by
$(a,0), (b,0)$ and $(0,1)$, $\Gamma \backslash G$ is a solvmanifold,
which is a finite quotient of $T^3$ and has a structure of $T^2$-bundle over
$T^1$.
\bigskip

\noindent {\bfseries Example 3.4.}
Let $\Lambda_n= \Z^2 \rtimes \Z$, where the action $\phi:
\Z \rightarrow {\rm Aut}(\Z^2)$ is defined by
$\phi(1)=A_n \in {\rm GL}(2,\Z)$,
$$A_n = \left(
\begin{array}{cc}
1 & n\\
0 & 1
\end{array}
\right).$$
The action $\phi$ can be extended to $\bar{\phi}: \R \rightarrow
{\rm Aut}(\R^2)$ defined by $\bar{\phi}(t)=
A(t) \in {\rm GL}(2,\R)$,
$$A(t) = \left(
\begin{array}{cc}
1 & t\\
0 & 1
\end{array}
\right).
$$
Let $N= \R^2 \rtimes \R$ be a nilpotent Lie group defined by
the action $\bar{\phi}$. Then each $\Lambda_n$ is a lattice of $N$,
and $\Lambda_n \backslash N$ is a nilmanifold.
\bigskip

\noindent {\bfseries Example 3.5.}
Let $\Gamma_n = \Lambda_n \rtimes \Z$, where $\Lambda_n$ is the
same nilpotent group as in Example 3.4, and the action $\phi:
\Z \rightarrow {\rm Aut}(\Lambda_n)$ is defined by $\phi(1)= \tau
\in {\rm Aut}(\Lambda_n)$,
$$\tau:  \left(
\begin{array}[c]{ccc}
1 & a & \frac{c}{n}\\
0 & 1 & b\\
0 & 0 & 1
\end{array}
\right) \;
\longrightarrow \;
\left(
\begin{array}[c]{ccc}
1 & a & \frac{-c}{n}\\
0 & 1 & -b\\
0 & 0 & 1
\end{array}
\right)
$$
Let $G= N \times \R$ be a nilpotent Lie group, where $N$ is the
same nilpotent Lie group as in Example 3.4. $\Gamma_n$  acts
as a group of automorphisms on $G$.
$\Gamma_n \backslash G$ is a solvmanifold, which has a nilmanifold with
$b_1 = 3$ as a double covering, and has a structure of $T^2$-bundle
over $T^2$.
\bigskip

\noindent{\large 3.2. \bfseries Complex Tori}

\smallskip

An $n$-dimensional torus $T^n$ is a compact homogeneous space of
the abelian Lie group $\R^n$: that is, $T^n = \Z^n \backslash \R^n$
where $\Z^n$ is an abelian lattice of $\R^n$ which is spanned by some
basis of $\R^n$ as a real vector space. For the case $n = 2m$,
the standard complex structure $\C^m$ on $\R^{2m}$ defines a complex
structure on $T^{2m}$. The complex manifold thus obtained is a {\em complex
torus}, which is known to be the only compact complex Lie group.
It should be noted that complex structures on $T^{2m}$ differ for different
abelian lattices, while the underlying differentiable manifolds are diffeomorphic
each other. On the other hand, as was shown in Example 3.2, $T^n (n \ge 3)$
can admits a structure of non-toral solvmanifold. It is conjectured that
complex structures on $T^{2m}$ are only the standard ones (we know that
this is valid for $m = 1,2$ since they are K\"ahlerian).
\bigskip

\noindent{\large  3.3. \bfseries Hyperelliptic surfaces}

\smallskip

Let $\Gamma = \Z^3 \rtimes \Z^1$, where the action
$\phi: \Z^1 \rightarrow {\rm Aut}(\Z^3)$ is defined by
$\phi(1) = A \in {\rm SL}(3,\Z)$. Assume that $A$ has a single root $1$,
and a double root $-1$ with linearly independent eigen vectors of $-1$
or non-real complex roots $\beta \,(|\beta|=1)$, $\bar{\beta}$.
We shall see that the class of the corresponding solvmanifolds with canonical
complex structures coincides with that of hyperelliptic surfaces.
\smallskip

For a given $A \in {\rm SL}(3,\Z)$ which satisfies our assumption,
we can find a basis $\{u_1, u_2, u_3\}$ of $\R^3$
such that $A u_1 = au_1 - bu_2, A u_2 = bu_1 + au_2, A u_3 = u_3$,
where $a = -1$ and $b = 0$ for the case that $A$ has a double root $-1$,
and $a ={\rm Re}\,\beta$ and $b = {\rm Im}\,\beta$ for the case that $A$ has a
non-real complex root $\beta$. 
Let $u_i=(u_{i1}, u_{i2}, u_{i3}),\; i=1,2,3$.
Then  $\{v_1, v_2, v_3 \},\;v_j=(u_{1j}, u_{2j}, u_{3j}), j=1,2,3$ defines an
abelian lattice $\Z^3$ of $\R^3$ which is preserved by a rotation around a fixed
axis. In particular, $\beta$ must be $e^{i \zeta}$
($\zeta = \frac{2}{3} \pi, \frac{1}{2} \pi$ or $\frac{1}{3} \pi$).

Furthermore, we may assume that $u_{3j} = 0, j = 1,2,$ and $A$ is of the form
$$\left(
\begin{array}[c]{ccc}
a_{11} & a_{12} & 0\\
a_{21} & a_{22} & 0\\
p & q & 1
\end{array}
\right),$$
where $A' = (a_{ij}) \in {\rm SL}(2, \Z), p, q \in \Z$.
Since $A'$ has  the root $-1$ (with linearly independent eigen vectors) or $\beta$,
we can assume that $A'$ is of the form: 
$$\left(
\begin{array}[c]{cc}
-1 & 0\\
0 & -1
\end{array}
\right),
\left(
\begin{array}[c]{cc}
0 & 1\\
-1 & -1
\end{array}
\right),
\left(
\begin{array}[c]{cc}
0 & 1\\
-1 & 0
\end{array}
\right),
\left(
\begin{array}[c]{cc}
0 & 1\\
-1 & 1
\end{array}
\right),
$$
according to the root $e^{i \eta}$ of $A'$, where
$\eta = \pi, \frac{2}{3} \pi, \frac{1}{2} \pi$ or $\frac{1}{3} \pi$, respectively. 

\smallskip

We now define a solvable Lie group $G = (\C \times \R) \rtimes \R$, where the
action $\phi: \R \rightarrow {\rm Aut}(\C \times \R)$ is defined by 
$$\phi(t)((z,s)) = (e^{i\eta t} z, s),$$
which is a canonical extension of $\phi$.
$\Gamma = \Z^3 \rtimes \Z$ clearly defines a lattice of $G$. 
Since the action on the second factor $\R$ is trivial, the multiplication of $G$ is
defined on $\C^2$ as follows:
$$(w_1,w_2) \cdot (z_1,z_2)= (w_1+e^{i \eta t} z_1, w_2+z_2),$$
where $t= {\rm Re}\,w_2$.  For each lattice $\Gamma$ of $G$,
$\Gamma \backslash G$ with canonical complex structure defines a complex surface,
actually a hyperelliptic surface. 
\smallskip

We can see that there exist seven isomorphism classes
of lattices of $G$, which correspond to seven classes of hyperelliptic surfaces.
For each $\eta$, take a lattice $\Z^3$ of $\R^3$ spanned by $\{v_1,v_2,v_3\}$ for
$A$ with $p = q = 0$. Then we can get a lattice $\Z^3$ spanned by $\{v_1,v_2,v_3'\}$
for $A$ with arbitrary $(p,q) \in \Z^2$, by changing $v_3$ into
$v_3' = s v_1 + t v_2 + v_3$ where $s, t \in \Q$ with $0 \le s,t < 1$, and
$\Gamma = \Z^3 \rtimes \Z$ defines a lattice of the solvable Lie group $G$.
By elementary calculation, we obtain the following seven isomorphism classes of
lattices: besides four trivial cases with $(p,q) = (0,0)$ and $(s,t) = (0,0)$ for
$\eta = \pi, \frac{2}{3} \pi, \frac{1}{2} \pi$ and $\frac{1}{3} \pi$,
we have three other cases with $(p,q) = (1,0)$, and 
(i) $(s,t) = (\frac{1}{2}, 0)$ for $\eta = \pi$,
(ii) $(s,t) = (\frac{1}{3}, \frac{1}{3})$ for $\eta = \frac{2}{3} \pi$,
(iii) $(s,t) = (\frac{1}{2}, \frac{1}{2})$ for $\eta = \frac{1}{2} \pi$.

The corresponding solvmanifolds with canonical complex structures
coincides with the class of all {\em hyperelliptic surfaces} (cf. [\ref{BPV}]).
\bigskip

\noindent{\large  3.4. \bfseries Primary Kodaira surfaces}

\smallskip

Let $\Gamma_n= \Lambda_n \times \Z$, and
$G= N \times \R$, where $\Lambda_n$ and $N$ are the
nilpotent group and nilpotent Lie group defined in Example 3.4
respectively. $S_n= \Gamma_n \backslash G$ is a nilmanifold
with $b_1=3$.  Expressing the nilpotent Lie group $N$ in Example 3.4 as
$$\left(
\begin{array}[c]{ccc}
1 & x & s\\
0 & 1 & y\\
0 & 0 & 1
\end{array}
\right),$$
we can define a coordinate change $\Phi$ from $N \times \R=
\R^3 \times \R$ to $\R^4$:
$$\Phi: ((x,y,s), t) \longrightarrow
(x,y, 2 s-x y, 2 t+ \frac{1}{2}(x^2+y^2)).$$
Considering $\R^4$ as $\C^2$ the group operation on $G$ in the new
coordinate can be expressed as follows:
$$(w_1,w_2) \cdot (z_1,z_2)=
(w_1 + z_1, w_2 - i \bar{w_1} z_1 + z_2).$$
$S_n$ with this complex structure is {\em a primary Kodaira surface}.

Conversely, any primary Kodaira surface, which is by definition a complex surface
with the trivial canonical bundle and $b_1=3$, can be written as $\C^2/\Gamma$
where $\Gamma$ is a properly discontinuous group of affine transformations
in the above form ([\ref{K2}]).
\bigskip

\noindent{\large 3.5.  \bfseries  Secondary Kodaira surfaces}

\smallskip

Let $\Gamma_n= \Lambda_n \rtimes \Z$, where the action $\phi:
\Z \rightarrow {\rm Aut}(\Lambda_n)$ satisfies the condition that
the induced automorphism $\widetilde{\phi}(1)$ of $\Z$ is trivial, that is,
$\widetilde{\phi}(1)= {\rm Id}$, and the induced automorphism
$\widehat{\phi}(1)$ of $\Z^2$ has a double root $-1$
with linearly independent eigen vectors, or non-real complex roots,
$\alpha, \bar{\alpha}\,
(\alpha \not= \bar{\alpha}),  |\alpha|= 1$.

We shall see that $\Gamma_n$ can be extended to a solvable Lie group
$G = N \rtimes \R$. As we have seen in Section 3.3, $\alpha$ must be
$e^{i \eta}$, $\eta = \pi, \frac{2}{3} \pi, \frac{1}{2} \pi$ or $\frac{1}{3} \pi$,
and there exists a basis $\{u_1', u_2'\}$ of $\R^2$ such that
$A u_1' = a u_1' - b u_2', A u_2' = b u_1' + a u_2'$, where $A = \widehat{\phi}(1)$,
$a = {\rm Re}\,\alpha, b = {\rm Im}\,\alpha$,
and $u_1' = (u_{11}, u_{12}), u_2' = (u_{21}, u_{22})$. 
The abelian lattice $\Z^2$ of $\R^2$ spanned by $\{v_1', v_2'\}$,
where $v_1' = (u_{11}, u_{21}), v_2' = (u_{12}, u_{22})$, is preserved by
the automorphism $\psi': (x, y) \rightarrow (ax - by, bx + ay)$ of $\R^2$.
We can extend $\psi'$ to an automorphism of $N$ of the form:
$$\psi: (x,y,z) \longrightarrow (a x-b y, b x+ a y, \frac{z}{n}+ h(x, y)),$$
for some polynomial $h$. As $\psi$ is a group homomorphism of finite order,
by simple calculation, we see that $h = \frac{1}{2} b\,(a x^2 - a y^2 - 2 b x y)$.
We can extend the lattice $\Z^2$ spanned by $\{v_1', v_2'\}$ to a lattice
$\Lambda_n$ spanned by $\{v_1, v_2, v_3\}, v_1 = (u_{11}, u_{21}, u_{31}),
v_2 = (u_{12}, u_{22}, u_{32}), v_3 = (0, 0, u_{33})$ for suitable $u_{31}, u_{32},
u_{33}$, so that $\Lambda_n$ is preserved by $\psi$. We now define a solvable
Lie group $G$ by extending the action $\phi(m) = \psi^m, m \in \Z$ to
$\psi(t), t \in \R$, replacing $a$ with $\cos \eta t$ and $b$ with $\sin \eta t$.
It is clear that $\Gamma_n = \Lambda_n \rtimes \Z$ defines a lattice of $G$.
If we take the new coordinate as in Section 3.4, the automorphism $\psi$ is
expressed as,
$$(z_1, z_2) \longrightarrow (\zeta z_1, z_2)$$
for $\zeta \in \C, |\zeta|=1$. It follows that the above automorphism
is holomorphic with respect to the complex structure defined in
Section 3.4. Therefore, $S_n= \Gamma_n \backslash G$ is a
finite quotient of a primary Kodaira surface, and $S_n$
with the above complex structure is {\em a secondary Kodaira surface}.

The classification of secondary Kodaira surfaces are known
(cf. [\ref{FM}, \ref{S2}]): it is by definition a finite quotient of 
a primary Kodaira surface, which can be also written as $\C^2/\Gamma$
where $\Gamma$ is a properly discontinuous group of affine transformations.
And all the secondary Kodaira surfaces are constructed in the above way. 
\bigskip

\noindent{\large  3.6. \bfseries Inoue surfaces of type $\bf S$}

\smallskip

Let $\Gamma= \Z^3 \rtimes \Z$, where the action $\phi:
\Z \rightarrow {\rm Aut}(\Z^3)$ is defined by $\phi(1)=
A \in {\rm Aut}(\Z^3)= {\rm GL}(3,\Z)$. Assume that $A$ has
complex roots $\alpha, \bar{\alpha}$ and a real root $c$,
where $c \not= 1, |\alpha|^2 c= 1$. Let $(\alpha_1,\alpha_2,\alpha_3)
\in \C^3$ be the eigen vector of $\alpha$ and  $(c_1,c_2,c_3)
\in \R^3$ the eigen vector of $c$. The set of vectors
$\{(\alpha_i,c_i) \in \C \times \R \,|\, i= 1,2,3\}$ are linearly
independent over $\R$, and defines a lattice $\Z^3$ of
$\C \times \R$. Let $G= (\C \times \R) \rtimes \R$ be a solvable
Lie group, where the action $\bar{\phi}: \R \rightarrow
{\rm Aut}(\C \times \R)$ is defined by
$\bar{\phi}(t): (z,s) \rightarrow (\alpha^t z,c^t s)$, which is
a canonical extension of $\phi$. Then $\Gamma$ is a lattice of $G$ and
$S= \Gamma \backslash G$ is a solvmanifold. Using the
diffeomorphism $\R \rightarrow \R_+$ defined by $t \rightarrow
e^{\log c t}$, $S$ can be considered as $\Gamma' \backslash
\C \times {\bf H}$, where $\Gamma'$ is
a group of automorphisms generated by $g_0$ and $g_i, i=1,2,3$,
which correspond to the canonical generators of $\Gamma$.
We can see that $g_0:(z_1,z_2) \rightarrow (\alpha z_1, c z_2)$ and
$g_i:(z_1,z_2) \rightarrow (z_1 + \alpha_i, z_2 + c_i), i=1,2,3$. $S$ with the
above complex structure is, by definition, {\em an Inoue surface of type $S$}.

\bigskip

\noindent{\large  3.7. \bfseries Inoue surfaces of type $\bf S^\pm$}

\smallskip

Let $\Gamma_n= \Lambda_n \rtimes \Z$, where the action $\phi:
\Z \rightarrow {\rm Aut}(\Lambda_n)$ satisfies the condition that
for the induced action $\widetilde{\phi}: \Z \rightarrow {\rm Aut}(\Z)$,
$\widetilde{\phi}(1)= {\rm Id}$, and for the induced action
$\widehat{\phi}: \Z \rightarrow {\rm Aut}(\Z^2)$,
$\widehat{\phi}(1)= (n_{ij}) \in {\rm SL}(2, \Z)$  has two positive real roots,
$a, b \;(ab=1)$. Let $(a_1,a_2), (b_1,b_2) \in \R^2$ be
eigen vectors of $a, b$ respectively. Let $G= N \rtimes \R$ be a solvable
Lie group, where the action $\bar{\phi}: \R \rightarrow {\rm Aut}(N)$
is defined by
$$\bar{\phi}(t): \left(
\begin{array}[c]{ccc}
1 & x & z\\
0 & 1 & y\\
0 & 0 & 1
\end{array}
\right)
\longrightarrow
\left(
\begin{array}[c]{ccc}
1 & a^t x & z\\
0 & 1 & b^t y\\
0 & 0 & 1
\end{array}
\right),
$$
which is a canonical extension of $\phi$.
In order to define a lattice $\Lambda_n$ which is preserved by
$\bar{\phi}$, we take $g_1,g_2,g_3 \in N$ as
$$g_1=\left(
\begin{array}[c]{ccc}
1 & a_1 & c_1\\
0 & 1 & b_1\\
0 & 0 & 1
\end{array}
\right),\;
g_2=\left(
\begin{array}[c]{ccc}
1 & a_2 & c_2\\
0 & 1 & b_2\\
0 & 0 & 1
\end{array}
\right),\;
g_3=\left(
\begin{array}[c]{ccc}
1 & 0 & c_3\\
0 & 1 & 0\\
0 & 0 & 1
\end{array}
\right),$$
where $c_1,c_2,c_3$ are to be determined, satisfying
the following conditions:

\begin{list}{}{\partopsep=0pt \parsep=0pt \leftmargin=20pt}
\item[1)]  $[g_1, g_2]=g_3^n$
\item[2)]  $\bar{\phi}(1)(g_1)=g_1^{n_{11}} g_2^{n_{12}} g_3^k,\;
\bar{\phi}(1)(g_2)=g_1^{n_{21}} g_2^{n_{22}} g_3^l$, where $k,l \in \Z.$
\end{list}
If we take $g_0 \in N \rtimes \R$ as
$$g_0=\left( \left(
\begin{array}[c]{ccc}
1 & 0 & p\\
0 & 1 & 0\\
0 & 0 & 1
\end{array}
\right), 1
\right),
\;p \in \R,$$
then $\{g_0,g_1,g_2,g_3\}$ defines a lattice $\Gamma_n$ of $G$, and
$S_n=\Gamma_n\backslash G$ is a solvmanifold.
\smallskip

Now, we define a diffeomorphism $\Phi: G=N \rtimes
\R \longrightarrow \R^3 \times \R_+$, for an arbitrary
$\gamma=p + qi \in \C$ and $\sigma=\log b$, by
$$\Phi: \left( \left(
\begin{array}[c]{ccc}
1 & y & x\\
0 & 1 & s\\
0 & 0 & 1
\end{array}
\right),
t \right)
\longrightarrow
(x,\, e^{\sigma t} y+ q \,t,\, s, \,e^{\sigma t}).$$
Then considering $\R^3 \times \R_+$ as $\C \times {\bf H}$,
$g_0,g_1,g_2,g_3$ are corresponding to the following holomorphic
automorphisms of $\C \times {\bf H}$,
$$g_0: (z_1, z_2) \rightarrow (z_1 + \gamma, b z_2),$$
$$g_i: (z_1, z_2) \rightarrow (z_1 + a_i z_2 + c_i, z_2 + b_i),$$
where $i=1,2, 3$ and $a_3 = b_3 = 0$.
$S_n$ with the above complex structure is, by definition,
{\em an Inoue surface of type $S^+$}.
\smallskip

An Inoue surface of type $S^-$ is defined similarly as
the case where the action $\phi: \Z \rightarrow {\rm Aut}(\Lambda_n)$
satisfies the condition that $\widetilde{\phi}(1)= -{\rm Id}$,
and $\widehat{\phi}(1)$ has a positive and a negative real root.
It is clear that an Inoue surface of type $S^-$ has $S^+$ with $\gamma=0$
as its double covering surface.
\bigskip

\noindent{\large 3.8.  \bfseries Complex structures on solvable Lie algebras}
\smallskip

We have seen that a four-dimensional solvmanifold $S$ can be expressed, up to
finite covering, as $\Gamma \backslash G$ where $G$ is a simply connected
solvable Lie group and $\Gamma$ is a lattice of $G$. And we have explicitly
constructed complex structures on $S$, as those canonically induced from some
left-invariant complex structures on $G$. In the following list, for each complex
surface with diffeomorphism type of solvmanifold, we express the corresponding
solvable Lie algebra $\mathfrak g$ of $G$ as being  generated by
$\{X_1,X_2,X_3,X_4\}$  with the specified bracket multiplication.  And for each case
except 6), the almost complex structure $J$ is defined by
$JX_1=X_2, JX_2=-X_1, JX_3=X_4, JX_4=-X_3$, for which the Nijenhuis tensor 
$N(X,Y)=[JX,JY]-J[JX,Y]-J[X,JY]-[X,Y]$ vanishes.
\smallskip

\begin{list}{}{\topsep=0pt \leftmargin=10pt \itemindent=5pt \parsep=0pt \itemsep=5pt}
\item[ 1)] Complex Tori \par
$[X_i, X_j]=0 \;( 1 \le i < j \le 4)$.
\item[ 2)] Hyperelliptic Surfaces \par
$[X_4, X_1]= -X_2, [X_4, X_2]= X_1$, and all other brackets vanish.
\item[ 3)] Primary Kodaira Surfaces \par
$[X_1, X_2]= -X_3$, and all other brackets vanish.
\item[ 4)] Secondary Kodaira Surfaces \par
$[X_1, X_2]= -X_3, [X_4, X_1]= -X_2, [X_4, X_2]= X_1$,
and all other brackets vanish.
\item[ 5)] Inoue Surfaces of Type $S$ \par
$[X_4, X_1]= a X_1 - b X_2, [X_4, X_2]= b X_1 + a X_2, [X_4, X_3]= -2a X_3$,
and all other brackets vanish, where $a, b \;(\not=0) \in \R$.
\item[ 6)] Inoue Surfaces of Type $S^+$ and $S^-$ \par
$[X_2, X_3]= -X_1, [X_4, X_2]= X_2, [X_4, X_3]= -X_3$, and all other
brackets vanish.
The almost complex structure $J$ is defined by $J X_1=X_2, J X_2=-X_1,
J X_3= X_4-q X_2, J X_4= -X_3-q X_1$, and the Nijenhuis tensor
vanishes for this $J$.
\smallskip

The Inoue Surfaces of Type $S^-$ is not of the form $\Gamma \backslash G$,
but has a double covering Inoue surface of type $S^+$; and its complex
structure comes from the one of the above type.
\end{list}

\bigskip

\begin{center}
{4.  PROOF OF THE MAIN THEOREM}
\end{center}

In this section we shall prove the main theorem we stated in Section 1;
we show actually that a complex surfaces with diffeomorphism type of
solvmanifold must be one of the following surfaces: complex tori,
hyperelliptic surfaces, Primary Kodaira surfaces, Secondary Kodaira surfaces,
and Inoue surfaces. It should be noted that we have already seen in section 3
that all of the above complex surfaces have the structures of four-dimensional
solvmanifold. In this section we shall use the standard notations
and terminologies in the field of complex surfaces (see [\ref{BPV}, \ref{FM}]
for references).
\smallskip

Now let $S$ be a complex surface with diffeomorphism type of
solvmanifold. We first remark that since $S$ is parallelizable (see Proposition 2.3)
the Euler number $c_2$ of $S$ vanishes, and the fundamental group of $S$
is abelian if and only if $S$ is a four-dimensional torus ([\ref{M}]).
Let $\kappa(S)$ be the Kodaira dimension of $S$. The classification of complex
surfaces with $c_2 = 0$ is divided into three cases: $\kappa(S) = -\infty, 0, 1$
(c.f. [\ref{BPV}]). In the case where $\kappa(S) = -\infty$, $S$ is a surface of
${\rm VII}_0$ or ruled surface of genus 1. The latter surface cannot be diffeomorphic
to a solvmanifold since the fundamental group of ruled surface of genus $1$
is $Z^2$ (which is abelian).
According to the well-known theorem of Bogomolov 
(proved by Li, Yau and Zheng [\ref{LYZ}]), we know that a complex surface of
${\rm VII_0}$ with $b_1 = 1$ and $c_2 = 0$ is an Inoue surface or a Hopf surface. 
Since the fundamental group of the Hopf surface is of the form $H \rtimes \Z$
where $H$ is a finite unitary group (including the trivial case) (cf. [\ref{H1}]),
it cannot be the fundamental group of solvmanifold. 
Hence $S$ must be an Inoue surface. In the case where $\kappa(S) = 0$,
$S$ is a complex torus, hyperelliptic surface, or Kodaira surface (primary or
secondary Kodaira surface).
\smallskip

In the case where $\kappa(S) = 1$, $S$ is, a (properly) elliptic surface (which is
minimal since $c_2 = 0$). Let us first recall some terminologies and fundamental
results concerning topology of elliptic surfaces in general. An elliptic surface is a
complex surface $S$ together with an elliptic fibration $f: S \rightarrow B$ where
$B$ is a curve, such that a general fiber $f^{-1}(t), t \in B$
(except finite points $t_1, t_2, ..., t_k$) is an elliptic curve. The base curve $B$
is regarded as a two-dimensional orbifold with multiple points $t_i$ with
multiplicity $m_i$, where $m_i \,(i \ge 2)$ is the multiplicity of the fiber
$f^{-1}(t_i)$ ($i = 1, 2, ..., k$). An elliptic surface $S$ is of the type {\em hyperbolic,
flat (Euclidian), spherical or bad}, according as the orbifold $B$ is of that type.
The Euler number $e^{orb}(B)$ of $B$ is by definition
$e(B) - \sum_{i = 1}^k 1 - \frac{1}{ m_i}$, where $e(B)$ is the Euler number of $B$
as topological space. We know that $B$ is hyperbolic, flat or spherical according
as $e^{orb}(B)$ is negative, 0, or positive. In the case $c_2 = 0$, we can
see ([\ref{W1}]) that an elliptic surface $S$ is hyperbolic, flat or spherical according
as the Kodaira dimesion $\kappa(S)$ is $1, 0$ or $-\infty$. 
We now continue our proof for the case $\kappa(S) = 1$. By the above argument,
$S$ is a minimal elliptic surface of hyperbolic type. We show that the
fundamental group of $S$ is not solvable, and thus $S$ cannot be diffeomorphic to
a solvmanifold. We have the following presentation of $\pi_1(S)$ as a short
exact sequence (c.f. [\ref{FM}]):
$$ 0 \rightarrow \Z^2 \rightarrow \pi_1(S) \rightarrow \pi_1^{orb}(B) \rightarrow 0,$$

\noindent where $\pi_1^{orb}(B)$ is the fundamental group as two-dimensional
orbifold.

Since $\pi_1^{orb}(B)$ is a discrete subgroup of $\rm PSL(2, \R)$, it contains
a torsion-free subgroup $\Gamma$ of finite index, such that $\Gamma$ is the
fundamental group of a finite orbifold covering $\widetilde{B}$ of $B$, which is a
close surface of genus $g \ge 2$ (or a Riemann surface of hyperbolic type)
(cf. [\ref{S}]). We know that $\Gamma$ is represented as a goup with generator
$\{a_1, a_2, ..., a_g, b_1, b_2, ..., b_g\}$ and relation $\prod_{i = 1}^g [a_i, b_i] = 1$,
which is not solvable for $g \ge 2$ (cf. [\ref{C}]).
It follows that $\pi_1(S)$ cannot be solvable since the quotient groups and
subgroups of a solvable group must be solvable. This completes the proof of
Main Theorem. 
\hfill $\Box$
\smallskip

We state the above result as the theorem, with its immediate consequence as
a proposition.

\bigskip

\noindent {\bfseries Theorem 4.1 (Main Theorem).}
{\em The complex surfaces with diffeomorphism type of solvmanifolds
are all of the complex tori, hyperelliptic surfaces,
Primary Kodaira surfaces, Secondary Kodaira surfaces and Inoue surfaces.}
\medskip

Combined with the results of section 3, we obtain the following result.
\smallskip

\noindent {\bfseries Proposition 4.2.}
{\em Every complex structure $J$ on a four-dimensional solvmanifold $S$ is
canonical: that is, up to finite covering, $S$ can be written as
$\Gamma \backslash G$, where $G$ is a simply connected solvable Lie group, $\Gamma$ is a lattice of $G$, and $J$ is the complex structure canonically induced
from a left-invariant complex structure on $G$}.
\medskip

\noindent {\bfseries Conjecture}.
It is natural to conjecture that the above proposition holds for solvmanifolds of
general dimensional. It should be noted that this conjecture is closely related to the
general conjecture on Kaehlerian solvmanifolds (see [\ref{H2}] for the detail).
\bigskip

\begin{center}
{5.  APPENDIX -- COMPLEX STRUCTURES ON FOUR-DIMENSIONAL\\
COMPACT HOMOGENEOUS SPACES}
\end{center}

In this section we determine complex surfaces with diffeomorphism type
of compact homogeneous spaces. 
It is known (due to V. V. Gorbatsevich [\ref{G}], also see [\ref{GO}]) that
a four-dimensional compact homogeneous space is diffeomorphic to one
of the following types: (1) $\prod S^{k_i}$ (up to finite quotient), where
$k_i \ge 1$ with $\sum {k_i} = 4$, (2) $\C {\rm P}^2$,
(3) Solvmanifold, (4) $S^1 \times \Gamma \backslash \widetilde{\rm SL}_2(\R)$,
where $\widetilde{\rm SL}_2(\R)$ is the universal covering of
${\rm SL}_2(\R)$ and $\Gamma$ is a lattice of $\widetilde{\rm SL}_2(\R)$.
We can determine, from the above result, complex structures on
four-dimensional compact homogeneous spaces in the following list:
\bigskip

\renewcommand{\multirowsetup}{\centering}

$$
\begin{array}{c|c|c|c}
\hline
S & b_1 & {\rm Complex\: Structure} & \kappa \\
\hline
S^2 \times T^2 & 2 & {\rm Ruled\: Surface\: of\: genus\: 1} &
\multirow{5}{7mm}{$\left. \rule{0mm}{10mm} \right\} {-\infty}$} \\ \cline{1-3}
S^1 \times_{\Z_m} S^3/H & 1 & {\rm Hopf\: Surface} & \\ \cline{1-3}
S^2 \times S^2 & 0 & {\rm Hirzebruch\: Surface\: of\: even\: type} & \\ \cline{1-3}
\C {\rm P}^2 & 0 & {\rm Complex\: Projective\: Space} & \\ \cline{1-3}
\multirow{5}{25mm}{\hbox{Solvmanifold $\left. \rule{0mm}{10mm} \right\{$}} & 1 &
{\rm Inoue\: Surface} & \\ \cline{2-4}
& 4 & {\rm Complex\: Torus} & 
\multirow{4}{7mm}{$\left. \rule{0mm}{8mm} \right\} 0$} \\ \cline{2-3}
& 3 & {\rm Primary\: Kodaira\: Surface} & \\ \cline{2-3}
& 2 & {\rm Hyperelliptic\: Surface} & \\ \cline{2-3}
& 1 & {\rm Secondary\: Kodaira\: Surface}&  \\ \hline
S^1 \times \Gamma \backslash \widetilde{\rm SL}_2(\R)& {\rm odd} &
{\rm Properly\: Elliptic\: Surface} & 1 \\
\hline

\end{array}
$$

\noindent where $\kappa$ is the Kodaira dimension of $S$, $\rm H$ is a finite
subgroup of ${\rm SU}(2)$ acting freely on $S^3$.
\medskip

It is well-known (due to A. Borel and J. P. Serre) that $S^4$ has no almost
complex structure. 
For the case of $S^2 \times T^2$, it is known
(due to T. Suwa [\ref{S1}]) that a complex surface is diffeomorphic to a
$S^2$-bundle over $T^2$ if and only if it is a ruled surface of genus $1$.
We can see this also from the recent result
(due to R. Friedman and Z. B. Qin [\ref{FQ}]) that the Kodaira dimension
of complex algebraic surface is invariant up to diffeomorphism.
To be more precise, there exist two diffeomorphism types of ruled
surfaces of genus $1$: the trivial one and the non-trivial one (which correspond
to two diffeomorphism types of $S^2$-bundles over $T^2$), and the
latter is not of homogeneous space form (see [\ref{S1}]).
For the case of $S^1 \times S^3$, K.~Kodaira showed [\ref{K1}]
that a complex surface diffeomorphic to a finite quotient of $S^1 \times S^3$ is a
Hopf surface. Generally a Hopf surface is diffeomorphic to a fiber bundle
over $S^1$ with fiber $S^3/{\rm U}$, defined by the action
$\rho: \pi_1(S^1) \rightarrow N_ {{\rm U}(2)}({\rm U})$
with $\rho(1)$ being cyclic of order $m$, where $\rm U$ is a finite subgroup of
${\rm U}(2)$ acting freely on $S^3$: that is, $S = S^1 \times_{\Z_m} S^3/U$
(cf. [\ref{H1}]).
We can see that a Hopf surface is of homogeneous space form if and only if
$U$ is a finite subgroup of ${\rm SU}(2)$. Let $G = {\rm SU}(2) \times S^1$,
which is a compact Lie group structure on $S^3 \times S^1$.
Take a finite subgroup $\Delta = H \rtimes \Z_m$ of $G$,
where $H$ is a finite subgroup of ${\rm SU}(2)$, $\Z_m$ is a finite
cyclic subgroup of $G$ generated by $c$:
$$c= (\tau, \xi),
\tau = \left(
\begin{array}{cc}
\xi^{-1} & 0\\
0 & \xi
\end{array}
\right), \xi^m = 1,$$
and $\tau$ belongs to $N_{{\rm SU}(2)}(H)$.
$S$ is a fiber bundle over $S^1$ with fiber $S^3/H$, which has a canonical
complex structure, defining a Hopf surface. It should be noted that if $\tau$
does not belong to $H$ and $m \ge2$, then $S$ is a non-trivial bundle.
Conversely, given a Hopf surface $S$ with fiber $S^3/H$, defined by the action
$\rho$, we can assume ([\ref{H1}]) that $\rho(1)$ is a diagonal matrix, all of which
entries are $m$-th roots of $1$. Then we can see that $\rho(1)$, which belongs
to $N_{{\rm U}(2)}(H)$, actually belong to $N_{ {\rm SU}(2)}(H)$.
Hence $S$ is diffeomorphic to the one constructed above.
For the case of $S^2 \times S^2$, it was shown
(due to Z. B. Qin [\ref{Q}]) that a complex surface diffeomorphic to
$S^2 \times S^2$ must be a Hirzebruch surface of even type, which is
by definition a ruled surface of genus $0$ with diffeomorphism type
$S^2 \times S^2$. As is well-known, there exist two  diffeomorphism types of
ruled surfaces of genus $0$: $S^2 \times S^2$ and $\C {\rm P}^2 \#
\overline{\C {\rm P}^2}$
(which correspond to two diffeomorphism types of $S^2$-bundles over $S^2$).
A Hirzebruch surface of odd type is the surface of the latter type. We can see that
no non-trivial finite quotient of $S^2 \times S^2$ has complex structure.
It is well-known (due to S. -T. Yau) that $\C {\rm P}^2$ can have only the
standard complex structure.
We have studied in detail the case of solvmanifolds in this paper.
The complex surfaces with diffeomorphism type of solvmanifolds are Inoue
surfaces for $\kappa = -\infty$ and all of those with $c_2=0$ for $\kappa = 0$.
For the case of $S^1 \times \Gamma \backslash \widetilde{\rm SL}_2(\R)$,
C.~T.~C.~Wall showed [\ref{W1}] that it admits a canonical complex
structure, which define a properly elliptic surface with $b_1={\rm odd}$
and $c_2=0$; and conversely any such surface with no singular fibers
is diffeomorphic to $S^1 \times \Gamma \backslash \widetilde{\rm SL}_2(\R)$
for some lattice $\Gamma$.

\bigskip
\noindent{\bfseries Acknowledgments}

The author is very grateful to S. Kobayashi and T. Ohsawa for their stimulating
suggestions and useful comments, and to Y. Matsushita for useful conversations
by e-mail. The author would like to thank V. V. Gorbatsevich for his very valuable
comments; in particular the results in Appendix have been obtained while
communicating with him by e-mail. Finally the author would like to express many
thanks to the referee for his careful reading of the earlier version of manuscript with
several important remarks, which lead to many improvements in the revised version.
\bigskip

\begin{center}
{REFERENCES}
\end{center}
\baselineskip=15pt
\renewcommand{\labelenumi}{[\theenumi]}
\begin{enumerate}
\itemsep=0pt
\item \label{A} L. Auslander, {\em An exposition of the structure of
solvmanifolds I II}, Bull. Amer. Math. Soc., {\bf 79} (1973),
No. 2, 227-261, 262-285.
\item \label{AS1} L. Auslander and R. H. Szczarba, {\em Characteristic
Classes of Compact Solvmanifolds}, Ann. of Math, {\bf 76} (1962), 1-8.
\item \label{AS2} L. Auslander and H. Szczarba, {\em Vector bundles
over noncompact solvmanifolds}, Amer. J. of Math., {\bf 97} (1975),
260-281.
\item \label{BPV} W. Barth, C. Peters and Van de Ven, {\em Compact Complex
Surfaces}, Ergebnisse der Mathematik und ihrer Grenzgebiete, Vol. 4,
Springer-Verlag, 1984.
\item \label{C} K. T. Chen, {\em Algeras of iterated path integrals and fundamental
groups}, Trans. Amer. Math. Soc., {bf 156}(1971), 359-379.
\item \label{FM} R. Friedman and J. W. Morgan, {\em Smooth Four-Manifolds
and Complex Surfaces}, Ergebnisse der Mathematik und ihrer Grenzgebiete,
Vol. 27, Springer-Verlag, 1994.
\item \label{FQ} R. Friedman and Z. B. Qin, {\em On complex surfaces
diffeomorphic to rational surfaces}, Invent. Math., {\bf 120} (1995), 81-117.
\item \label{LYZ} J. Li, S. T. Yau and F. Zheng, {\em A simple proof of
Bogomolov's theorem on class $VII_0$ surfaces with $b_2=0$},
Illinois J. Math., {\bf 34} (1990), 217-220.
\item \label{G} V. V. Gorbatsevich, {\em On the classification of
four-dimensional compact~homogeneous spaces}, Usp. Mat. Nauk.,
{\bf 32} (1977), No. 2, 207-208 (Russian).
\item \label{GO} V. V. Gorvatsevich and A. L. Onishchik, {\em Lie Transformation
Groups}, Encyclopaedia of Mathematical Sciences, Vol. 20, Lie Groups and
Lie Algebras I, Springer-Verlag, 1994.
\item \label{H1} K. Hasegawa, {\em Deformations and diffeomorphism types of
generalized Hopf manifolds}, Illinois J. Math., {\bf 37} (1993), 643-651.
\item \label{H2} K. Hasegawa, {\em A class of compact K\"ahlerian
solvmanifolds and a general conjecture}, Geometriae Dedicata,
{\bf 78} (1999), 253-258.
\item \label{I} M. Inoue, {\em On surfaces of class $VII_0$},
Inventiones math., {\bf 24} (1974), 269-310.
\item \label{K1} K. Kodaira, {\em Complex structures on $S^3 \times S^1$},
Proc. Nat. Acad. Sci., {\bf 55} (1966), 240-243.
\item \label{K2} K. Kodaira, {\em On the structure of compact complex analytic
surfaces, I.}, Amer. J. Math., {\bf 86} (1964), 751-798; {\em II}., Am. J. Math.,
{\bf 88} (1966), 682-721.   
\item \label{M} G. D. Mostow, {\em Factor spaces of solvable Lie groups},
Ann. of Math., {\bf 60} (1954), 1-27.
\item \label{Q} Z. B. Qin, {\em Complex structure on certain differential
$4$-manifolds}, Topology, {\bf 32} (1993), No. 3, 551-566.
\item \label{S} P. Scott, {\em The Geometries of 3-manifolds},
Bull. London Math. Soc., {\bf 15} (1983), 401-487.
\item \label{SF} K. Sakamoto and S. Fukuhara, {\em Classification of
$T^2$ bundles over $T^2$}, Tokyo J. Math., {bf 6} (1983), 311-327.
\item \label{S1} T. Suwa, {\em Ruled surfaces of genus $1$},
J. Math. Soc. Japan, {\bf 21} (1969), 291-311.
\item \label{S2} T. Suwa, {\em Compact quotient spaces of $\C^2$ by affine
transformation groups}, J. Differential Geometry, {\bf 10} (1975), 239-252. 
\item \label{U} M. Ue, {\em Geometric 4-manifolds in the sense of
Thurston and Seifert 4-manifolds I}, J. Math. Soc. Japan,
{\bf 42} (1990), 511-540.
\item \label{W1} C. T. C. Wall, {\em Geometric structures on compact complex
analytic surfaces}, Topology, {\bf 25} (1986), No. 2, 119-153.
\item \label{W2} C. T. C. Wall, {\em Geometries and geometric structures
in real dimension $4$ and complex dimension $2$}, Lecture Notes
in Math., {\bf 1167} (1986), 118-153.
\item \label{WG} H. C. Wang, {\em Discrete subgroups of solvable Lie
groups I}, Ann. of Math., {\bf 64} (1956), 1-19.
\end{enumerate}
\bigskip
\baselineskip12pt
\begin{flushleft}
Keizo Hasegawa\\
\vskip1pt
Department of Mathematics\\
Faculty of Education and Human Sciences\\
Niigata University, Niigata\\
JAPAN\\
\vskip3pt
e-mail: hasegawa@ed.niigata-u.ac.jp\\
\end{flushleft}
\end{document}